\documentclass{article}
\usepackage{amssymb}
\usepackage{amsmath}
\usepackage{amsfonts, latexsym}

\begin{document}

\title{A Koksma-Hlawka-Potential Identity on the $d$ Dimensional Sphere and its Applications to Discrepancy}
\author{Steven B. Damelin \thanks{Mathematical Reviews, 416 Fourth Street, Ann Arbor, MI 48103, USA\, email: damelin@umich.edu} }
\date{}
\maketitle

\begin{abstract}
Let $d\geq 2$ be an integer, $S^d\subset {\mathbb R}^{d+1}$ the unit sphere and
$\sigma$ a finite signed measure
whose positive and negative
parts are supported on $S^d$ with finite energy.  In this paper, we derive an error estimate for
the quantity $\left|\int_{S^d}fd\sigma\right|$, for a class of
harmonic functions $f:\mathbb R^{d+1}\to \mathbb R$. Our error estimate 
involves
2 sided bounds for a Newtonian potential with respect to $\sigma$ away from 
its support.
In particular, our main result allows us to study quadrature errors, for
scatterings on the
sphere with given mesh norm.
\end{abstract}
\bigskip



\setcounter{equation}{0}
\section{Introduction}

Integration and discrepancy are important problems in applied mathematics and approximation theory
and in many applications, one needs to estimate the quantity
${\rm sup}_{f\in {\cal F}}\left|\int_{{\cal B}} fd\zeta\right|$
where ${\cal B}\subset \mathbb R^{d+1}$ is a bounded domain or manifold, $d\geq 2$
is an integer, $\zeta$ is a Borel measure
with compact support in ${\cal B}$ and ${\cal F}$ is a suitable class of
real valued functions with domain ${\cal B}$. Such problems above, arise naturally in many
interdisciplinary areas such as mathematical finance, imaging, geodesy,
scattering and statistical learning theory.

Our main objective in this paper is as follows. We derive an error estimate for
the quantity $\left|\int_{S^d}fd\sigma\right|$, for a class of
harmonic functions $f:\mathbb R^{d+1}\to \mathbb R$. Here and henceforth,
$S^d$ will denote the $d$ dimensional sphere realized as a subset of $\mathbb R^{d+1}$ and
$\sigma$ will denote a finite signed measure,
whose positive and negative parts are supported on $S^d$ with finite energy.
Our error estimate involves 2 sided bounds for a Newtonian potential with
respect to $\sigma$ away
from its support. In particular, our main result allows us to study quadrature
errors for scatterings on the
sphere with given mesh norm.
\bigskip

\setcounter{equation}{0}
\section{Spherical harmonics, multipole expansions, potentials and pseudo-differential operators}

In this section, we introduce needed notation and
pertinent facts concerning
spherical harmonics, multipole expansions, potentials and pseudo-differential
operators which we use in the
sequel. We also state our main results.

\subsection{Spherical harmonics}

In this subsection, we collect together some pertinent facts re spherical
harmonics which we will need throughout.
\medskip

Let $d\geq 2$ be given and define
\[
S^d:=\left\{x:=(x_1,...,x_{d+1})\in \mathbb R^{d+1}:\, x_1^2+...+x_{d+1}^2=1\right\}
\]
to be the surface of the $d$ dimensional unit sphere realized as a subset of $\mathbb R^{d+1}$
Euclidean space. By $rS^d$ we shall mean the surface of a sphere
of radius $0<r_0<r\leq 1$ for some fixed positive $r_0$. By $B^d(0,\delta)$, we will always mean the
open ball of radius
$0<\delta<1$ in $\mathbb R^{d+1}$ and by $\overline{B^d(0,\delta)}$ its closure with
boundary $\delta S^d$. Once and for all, let
$\sigma_d:=\frac{\mu_d|_{S^d}}{w_d}$,
denote normalized Lebesgue measure in $\mathbb R^{d+1}$ restricted to
$S^d$ where
$w_d$ denotes the volume of $S^d$ given by
\[
\frac{2\pi^{\frac{d+1}{2}}}{\Gamma(\frac{d+1}{2})}
\]
where $\Gamma$ is the gamma function.  Henceforth, ${\cal M}_d$, will
denote the class of all finite signed
measures $\sigma$ whose positive and negative parts are supported
on $S^d$ with finite energy. Given $f:\mathbb R^{d+1}\to \mathbb R$ and any
$0<r_0<r\leq 1$, we have an associated norm
\[
||f||_{L_p(rS^d)}:=
\left\{ {
\begin{array}{ll}
\left(\int_{rS^d}|f(x)|^pd\mu_d(x)\right)^{1/p}, &
1\leq p<\infty \\
{\rm esssup}_{x\in rS^d}|f(x)|, & p=\infty.
\end{array} }
\right.
\]
The class of all measurable functions $f:rS^d\to \mathbb R$ for which
$||f||_{L_{p}(rS^d)}<\infty$
will be denoted by $L_{p}(rS^d)$, with the usual
understanding that functions that are equal
almost everywhere are considered equal elements of $L_{p}(rS^d)$.
\medskip

The usual inner product of 2 vectors $x,y\in \mathbb R^{d+1}$ will be denoted
by $x.y$ so that
\[
S^d:=\left\{x\in \mathbb R^{d+1}:\, ||x||=1\right\}
\]
where $||x||=(x.x)^{1/2}$ is the Euclidean norm of $x\in \mathbb R^{d+1}$.
and for any $x,y\in S^d$, $x.y\in [-1,1]$.
\medskip

For a fixed integer $l\geq 0$, the restriction to $S^d$ of a
homogeneous harmonic polynomial of degree $l$, is called a
spherical harmonic of degree $l$.
The dimension of the space of spherical harmonics of degree $l$,
which we denote by $Z(d,l)$, is given by:
\begin{equation}\label{2.1}
Z(d,l):=
\left\{ {
\begin{array}{ll}
\frac{2l+d-1}{l+d-1}{l+d-1
\choose l}, & l\geq 1 \\
1, & l=0.
\end{array} }
\right.
\end{equation}
Let
\[
\left\{Y_{l,k}:l=0,1,...;\,k=1,...,Z(d,l)\right\}
\]
be a real orthonormal basis for $L^2(S^d)$.
\medskip

We have for each $x,y\in S^d$, the well known
\medskip

{\bf Addition formula}:
\begin{equation}\label{eq 2.2}
\sum_{k=1}^{Z(d,l)}Y_{l,k}(x)Y_{l,k}(y)=
\frac{Z(d,l)}{\omega_d}P_l(d+1,x.y), l=0,1,....,
\end{equation}
where $P_l(d+1,.)$ is the Legendre polynomial of degree $l$ in $d+1$ dimensions
over $[-1,1]$.

The Legendre polynomials are normalized so that $P_l(d+1,1)=1$ for
each $l$ and they satisfy the
orthogonality relations:
\[
\int_{-1}^1P_l(d+1,x)P_s(d+1,x)(1-x^2)^{d/2-1}dx=
\frac{\omega_d \delta_{l,s}}
{\omega_{d-1}Z(d,l)}.
\]
We will need their relations to the Gegenbauer polynomials
$P_l^{(d-1)/2}$ given by
\begin{equation}\label{eq 2.3}
P_l^{(d-1)/2}(x)={l+d-2\choose l}P_l(d+1,x),\, l\geq 0,\, x\in [-1,1]
\end{equation}
and the useful fact that
\begin{equation}\label{eq 2.4}
(1-2rt+r^2)^{-(d-1)/2}=\sum_{l=0}^{\infty}r^lP_{l}^{(d-1)/2}(t),\, 0<r_0<r<1,\, |t|\leq 1.
\end{equation}
In particular,  if $x\in rS^d$ for some $r_0<r<1$ and $\eta \in S^d$, then it is well known that we may write
$x=r\zeta,\, \zeta\in S^d$ and obtain the formula:
\begin{equation}\label{eq 2.5}
\frac{1}{||x-\eta||^{d-1}}= \frac{1}{(1+r^2-2r(\zeta.\eta))^{(d-1)/2}}
=\sum_{l=0}^{\infty}r^lP_{l}^{(d-1)/2}(\zeta.\eta).
\end{equation}
\medskip

In the sequel, we will also make use of the following:
\medskip

{\bf Funk-Hecke Formula}:\, Let $f:[-1,1]\to\mathbb R$ be a continuous function. Then
\begin{equation}\label{eq 2.6}
\int_{S^d}f((\eta.\zeta))Y_{l,k}(\zeta)d\sigma_d(\zeta)=\lambda Y_{l,k}(\eta),\,
\eta\in S^d
\end{equation}
where $l=0,1,...;\,k=1,...,Z(d,l)$ and
\[
\lambda:=\frac{\omega_{d-1}}{\omega_{d}}\int_{-1}^1 f(t)P_l(d+1,t)(1-t^2)^{(d-2)/2}dt.
\]

We adopt the following convention.
Given $f:\mathbb R^{d+1}\to \mathbb R$ and $0<r_0<r\leq 1$,
we let $f_r:S^{d}\to \mathbb R$ be defined by
\[
f_r(x):=f(rx),x\in S^d.
\]
Then we write
\begin{equation}\label{eq 2.7}
\hat{f_r}(l,k):=\int_{S^d}f_r(x)Y_{l,k}(x)d\mu_d(x), \, l\geq 0, 1\leq k\leq Z(d,l)
\end{equation}
where
\begin{equation}\label{eq 2.8}
\sum_{l=0}^{\infty}\sum_{k=1}^{Z(d,l)}\hat{f_r}(l,k)Y_{l,k}(x),\, x\in S^{d}
\end{equation}
 is the spherical expansion of $f_r$ on $S^d$.

\subsection{Measure of Discrepancy}

As our measure of discrepancy, we will use
the Newtonian potential with respect to the class  of measures $\sigma\in {\cal M}_d$ given by
\begin{equation}\label{eq 2.9}
U^{\sigma}(x):=\int_{S^d}\frac{1}{||x-y||^{d-1}}d\sigma(y), x\in \mathbb R^{d+1}.
\end{equation}
Henceforth, (we recall) we write
\[
\sigma:=\sigma^{+}-\sigma^{-}
\]
where both measures $\sigma^{\pm}$ are finite, non negative and supported on $S^d$ and with finite energy.
It is well known that $U^{\sigma^{\pm}}$ exist, are locally
integrable and finite almost
everywhere with respect to $d$ dimensional Hausdorf measure. Moreover,
$U^{\sigma^{\pm}}$
are superharmonic in $\mathbb R^{d+1}$ and harmonic outside the support of
$\sigma^{\pm}$, see \cite{L}.
We will need to estimate $U^{\sigma}$ away from the support of $\sigma$ and to
this end, we will need to take the balayage
of $\sigma^{\pm}$ onto a suitably defined region in $\mathbb R^{d+1}$.
More precisely, given $0<r_0<r<1$, we set
\begin{equation}\label{eq 2.10}
G:=\mathbb R^{d+1}\backslash \overline{B^d(0,r)}\cup \{\infty\}
\end{equation}
with compact boundary $\partial G=rS^d$. Observe that the supports of
$\sigma^{\pm}$ are contained in $\overline{G}$.
Let $\sigma_B^{+}$ denote the balayage of $\sigma^{+}$ onto $\partial G$
and let $\sigma_B^{-}$ denote the balayage of $\sigma^{-}$ onto $\partial G$.
As $G$ is regular with respect to the Dirichlet problem on $\mathbb R^{d+1}$, see
\cite[Theorems 4.2, 4.5]{L}, $\sigma_B^{\pm}$ exist,
are unique and have the following two additional properties:
\begin{equation}\label{eq 2.11}
U^{\sigma^{\pm}}(x)=U^{\sigma_B^{\pm}}(x),\, x\in rS^d
\end{equation}
and
\begin{equation}\label{eq 2.12}
\int_{S^d}f d\sigma^{\pm}=
\int_{\overline{G}}fd\sigma^{\pm}=\int_{\partial G}fd\sigma_{B}^{\pm}=
\int_{rS^d}fd\sigma_B^{\pm}
\end{equation}
for all functions $f:\mathbb R^{d+1}\to\mathbb R$,
harmonic on $G$, continuous on ${\overline{G}}$ and satisfying
$f(\infty)=0$.
\medskip

Here and throughout, we will adopt the convention that $C,C_1,...$ will
always denote positive, finite
constants independent of $f$, $x$, $r$ and $n$ but possibly depending
on other parameters such as $d$, $\varepsilon$, $s$, $||\sigma||$, $\alpha$ and $r_0$.
These constants may also take on different values at different times.
Throughout, lower case summation indices such as $l$, $k$, $m$ and $j$ will run through
a subset of nonnegative integers unless stated otherwise.

\subsection{Space of Functions}

In this subsection, we define our approximation space.
\medskip

Given $s>\frac{d}{2}$ and $0<r_0<r\leq 1$, we shall henceforth say that
$f_r\in L_1(S^d)$ is an element of the space $H_s(S^d)$ if
\begin{equation}\label{eq 2.13}
\sum_{l=0}^{\infty}\sum_{k=1}^{Z(d,l)}(\hat{f_r}(l,k))^2m_l^2<\infty
\end{equation}
where
\begin{equation}\label{eq 2.14}
m_l=m(l,d,s):=
\left\{ {
\begin{array}{ll}
\frac{l^s(2l+d-1)}{(d-1)\omega_d}, & l\geq 1 \\
1, & l=0
\end{array} }
\right.
\end{equation}
We may define a norm on the space $H_s(S^d)$ by setting:
\begin{equation}\label{eq 2.15}
||f_r||_{H_s(S^d)}:= \left(\sum_{l=0}^{\infty}
\sum_{k=1}^{Z(d,l)}(\hat{f_r}(l,k))^2m_l^2\right)^{1/2}.
\end{equation}
\medskip

{\bf D-Operator}\, A fundamental tool in our analysis will be the following operator which
we now introduce. Given $f_r\in L_1(S^d),\, 0<r_0<r\leq 1$, we define (formally) the operator
\begin{equation}\label{eq 2.16}
D(f_r):=\sum_{l=0}^{\infty}\sum_{k=1}^{Z(d,l)}
\frac{(2l+d-1)}{(d-1)\omega_d}\hat{f_r}(l,k)Y_{l,k}.
\end{equation}
\bigskip

We begin with our first basic result concerning
the class $H_s(S^d)$ and the operator $D$ given by (2.16).
\medskip

{\bf Proposition 1}\, Let $d\geq 2$, $0<p\leq \infty$, $0<r_0<r\leq 1$ and $s>\frac{d}{2}$. Then the
following hold true:
\begin{itemize}
\item[(i)] Let $f_r\in H_s(S^d)$ and suppose $f_r$ can be recovered pointwise
by its spherical expansion on $S^d$. Then
\begin{equation}\label{eq 2.17}
||f_r||_{L_{p}(S^d)}\leq C^{*}||f_r||_{H_s(S^d)}
\end{equation}
where
\[
C^{*}:=\left[\sum_{l=1}^{\infty}\frac{e^d\omega_{d}(d-1)^2l^{d-1-2s}}{(2l+d-1)^2}+
\frac{1}{\omega_d}\right]^{1/2}.
\]
\item[(ii)] Uniformly for $f_r\in H_s(S^d)$, we have
\begin{equation}\label{eq 2.18}
||D(f_r)||_{L_{p}(S^d)}\leq C^{**}||f_r||_{H_s(S^d)}
\end{equation}
where
\[
C^{**}:=\left[\sum_{l=1}^{\infty}\frac{e^dl^{d-1-2s}}{\omega_d}+\frac{1}{\omega_d^3}
\right]^{1/2}\omega_d.
\]
\item[(iii)] Suppose that $s>\frac{3d-2}{4}$ and let $f_r\in H_s(S^d)$. Suppose $f_r$ can be recovered
pointwise by its spherical expansion on $S^d$. Then $f_r$ is Lipschitz of order $1$ with
Lipschitz constant
\begin{equation}\label{eq 2.19}
C||f_r||_{H_s(S^d)}
\end{equation}
for some explicit positive constant $C$ depending on $d$ and $s$.
\end{itemize}
Henceforth, given fixed $d\geq 2$ and $0<r_0<r<1$, ${\cal F}_r$ will
denote the class of functions satisfying the following:
\begin{itemize}
\item $f:\mathbb R^{d+1}\to \mathbb R$ is harmonic on $G$, continuous
on $\overline{G}$ and vanishes at infinity.
\item
\[
||D(f)||_{L_{\infty}(rS^d)}<\infty.
\]
\end{itemize}
\medskip

\subsection{Discrepancy Results}

Our main discrepancy result is given in:
\medskip

{\bf Theorem 2:}\, Let $d\geq 2$, $0<r_0<r<1$ and choose $f\in {\cal F}_r$.
Also let
$1\leq p,p'\leq \infty$ with $\frac{1}{p}+\frac{1}{p'}=1$ and
$\sigma\in {\cal M}_d$. Suppose, moreover that
\[
f_r(\eta)=\sum_{l,k}\hat{f}_r(l,k)Y_{l,k}(\eta), \, \eta\in
S^d.
\]
Then we have:
\begin{eqnarray}\label{eq 2.20}
&& \left|\int_{S^d}f(x)d\sigma(x)\right|\leq \\
&& \nonumber \leq \frac{1}{r_0}||D(f)||_{L_{p}(rS^d)}
||U^{\sigma}||_{L_{p'}(rS^d)}.
\end{eqnarray}
\medskip

{\bf Remark}
\begin{itemize}
\item[(a)] Our estimate consists of 2 contributions. The second contribution
is a discrepancy estimate in terms of 2 sided bounds for a
Newtonian potential with respect to $\sigma$ on an inner sphere away
from the support of $\sigma$. As we show below, this later quantity may be estimated
for certain measures defining point systems with given mesh norm. The first
contribution, involves an $L_p$ norm of $D(f)$ on the sphere $rS^d$.
\end{itemize}
\medskip
If we now specialize the measures in Theorem 2, we obtain:
\medskip

{\bf Corollary 3:}\, Let $d\geq 2$, $0<r_0<r<1$ and choose $f\in {\cal F}_r$
satisfying the condition of Theorem 2.
Let $E_o$ be a scattering of $n\geq 1$ distinct points $t_{k,n},\, 1\leq k\leq n$
on $S^d$ and $a_{k,n},\, 1\leq k\leq n$, $n$ real weights. Define
\begin{equation}\label{eq 2.21}
\nu_n(x):=\sum_{k=1}^na_{k,n}\delta_{t_{k,n}}(x),\, x\in S^d
\end{equation}
where $\delta_{t_k}(.)$ denotes Dirac mass and let $\mu$ be a Borel measure
on $\mathbb R^{d+1}$ with
support in $S^d$. Let $1\leq p,p'\leq \infty$ with $\frac{1}{p}+\frac{1}{p'}=1$. Then
\begin{eqnarray}\label{eq 2.22}
&& \left|\int_{S^d}fd\mu-\sum_{k=1}^na_{k,n}f(t_{k,n})\right| \\
\nonumber && \leq \frac{1}{r_0}||D(f)||_{L_{p}(rS^d)}
||U^{\mu-\nu_n}||_{L_{p'}(rS^d)}.
\end{eqnarray}
\medskip

We now focus on estimating the discrepancy term $||U^{\sigma}||_{L_{p'}(rS^d)}$
in Theorem 2 for natural choices of $\sigma$. For simplicity we will consider the
case $p'=\infty$.
Let us recall quickly that given any
finite scattering $E_o$ of
distinct points on $S^d$, the mesh norm of $E_o$ is defined by
\begin{equation}\label{eq 2.23}
\delta_{E_o}:={\rm max}_{x\in S^d}{\rm dist}(x,E_o).
\end{equation}
Moreover, if ${\cal R}$ denotes a finite partitioning of $S^d$, then the partition norm
for ${\cal R}$ is defined by
\begin{equation}\label{eq 2.24}
||{\cal R}||:={\rm max}_{R\in {\cal R}}({\rm diam}R).
\end{equation}
We have:
\medskip

{\bf Theorem 4}
\begin{itemize}
\item[(a)] Let $d\geq 2$, $0<r_0<r<1$ and $E_o$ a scattering of $n\geq 1$ distinct
points $t_{k,n},\, 1\leq k\leq n$
on $S^d$. Suppose there exists a finite disjoint partitioning ${\cal R}$ of the sphere into
$n$ subsets $R_{k,n},\, 1\leq k\leq n$ and with
$S^q=\cup_{k=1}^n R_{k,n}$ with each $R_{k,n}$ containing exactly one $t_{k,n}$.
Let $\mu$ be a finite Borel measure on $\mathbb R^{d+1}$ with support in $S^d$
and suppose $\nu_n$ is given by (2.22) with
$a_{k,n}=\mu(R_{k,n}),\, k=1,...,n$. Then setting $\sigma=\mu-\nu_n$, gives
\begin{equation}\label{eq 2.25}
||U^{\sigma}||_{L_{\infty}(rS^d)}\leq \frac{(d-1)||\mu||}{(1-r)^{d+1}}||{\cal R}||.
\end{equation}
\item[(b)] Let $\epsilon>0$, $d\geq 2$, $E_o$ a finite scattering of
$n\geq 1$ distinct points on $S^d$ and let $\mu$ be a finite Borel measure on
$\mathbb R^{d+1}$ with support in $S^d$. Suppose $E_o$ is chosen so that $\delta_{E_o}$ is
small enough to ensure that
\begin{equation}\label{eq 2.26}
\left(\frac{(d-1)8d\sqrt{2d(d+1)}||\mu||\delta_{E_o}}{\epsilon}\right)<1.
\end{equation}
Then there exists a reduction $E$ of $E_o$ and an explicit
finite disjoint partitioning
${\cal R}$ of $S^d$ with each $R\in {\cal R}$ containing one point of
$E$ and at least one point of $E_o$ so that
\begin{equation}\label{eq 2.27}
\left(\frac{(d-1)||\mu||||R||}{\epsilon}\right)<1.
\end{equation}
Moreover, if $\mu_{{\rm card}E}$ is defined as in (2.21)
with respect to $E$ and with positive weights
\[
a_{k,{\rm card}E},\, 1\leq k\leq {\rm card}E
\]
and in addition, we define the signed measure
$\sigma=\mu-\mu_{{\rm card}E}$, then the following holds:\, For
$0<r_0<r<1-\left(\frac{(d-1)||\mu||||R||}{\epsilon}\right)
^{\frac{1}{d+1}}$, we have
\begin{equation}\label{eq 2.28}
||U^{\sigma}||_{L_{\infty}(rS^d)}\leq \epsilon.
\end{equation}
\end{itemize}

{\bf Remark:}\,
\begin{itemize}
\item[(a)] Notice that if the mesh norm of $E_o$ is small,
(ie if for example the scattering $E_o$ is not concentrated at
one place on the sphere such as one of the poles), then we are able to control the estimate for the potential
in Theorem 2 with the help of Theorem 4.
\item[(b)] We mention that it is possible to show that, there exists a
disjoint (equal area) partitioning ${\cal R}$ of the sphere $S^d$ into $n$ parts such that
each $R_{k,n}\in {\cal R},\, 1\leq k\leq n$
satisfies
\begin{equation}\label{eq 2.29}
\mu_d(R_{k,n})\sim \frac{1}{|R_{k,n}|}=1/n,\, ||R_{k,n}||\leq \frac{C}{n^{1/d}}
\end{equation}
for some $C>0$ independent of $n$ but depending on $d$.
Thus provided we have $n$ points with exactly one point in each partition,
we may apply Theorem 4 to equal weighted quadrature rules as well.
\item[(c)] Requiring information on mesh norm is natural since there are many
interesting examples of point systems on the sphere where good estimates on mesh norm
are available, for example extremal fundamental systems.
\end{itemize}

The remainder of this paper is devoted to the proofs of
Theorems 2-4.

\setcounter{equation}{0}
\section{Proofs}

In this section, we prove our results. We begin with the proof of Theorem 2.
We find it instructive to briefly summarize the ideas of the proof.
Let $G$ be the region defined by (2.10).  As $\sigma\in {\cal M}_d$,
$\sigma^{\pm}$ are supported
on $S^d\subset \overline{G}$. Moreover, the balayage of
$\sigma^{\pm}$, $\sigma_B^{\pm}$ exist, are supported on
$rS^d=\partial G$ and satisfy
\[
\int_{S^d}fd\sigma^{\pm}=\int_{\overline{G}}fd\sigma^{\pm}=
\int_{\partial G}fd\sigma_{B}^{\pm}=\int_{rS^d}fd\sigma_B^{\pm}.
\]
Note that for the above string of equations to hold, it is enough
that $f$ be continuous on $\partial G$, harmonic on $G$ and take the value $0$ at infinity.
Next, recalling the operator $D$ defined by (2.16), we prove that
\[
\int_{rS^d}fd\sigma_B^{\pm}=\frac{1}{r}\int_{rS^d}D(f)U^{\sigma_B^{\pm}}d\mu_d
\]
so that using the property of balayage again, we have
\[
\int_{S^d}fd\sigma^{\pm}=\frac{1}{r}\int_{rS^d}D(f)U^{\sigma^{\pm}}d\mu_d
\]
Adding, we obtain
\[
\int_{S^d}fd\sigma=\frac{1}{r}\int_{rS^d}D(f)U^{\sigma}d\mu_d.
\]
\medskip

Theorem 2 will follow from the following lemma which is of independent interest.
\medskip

{\bf Lemma 1}\, Let $d\geq 2$, $\sigma\in {\cal M}_d$, $0<r_0<r<1$ and assume
$f$ satisfies the assumptions of Theorem 2.
Then the following holds.
\begin{eqnarray}\label{eq 3.1}
&& \int_{S^d}r^{1-d}fd\sigma^{\pm}=
\int_{rS^d}r^{-d}D(f)U^{\sigma_B^{\pm}}d\mu_d\\
\nonumber && =\int_{rS^d}r^{-d}D(f)U^{\sigma^{\pm}}d\mu_d.
\end{eqnarray}
\medskip

{\bf Proof}\, For notational convenience, we shall henceforth write
\[
\sum_{l=0}^{\infty}\sum_{k=1}^{Z(d,l)}=\sum_{l,k}.
\]
The first thing to do, is to calculate the coefficients of
the expansion
of $\sigma_B^{\pm}$. Thus to see this, lets recall, see (2.11), that we have
\[
\int_{rS^d}\frac{d\sigma_B^{\pm}(y)}{||x-y||^{d-1}}=
\int_{S^d}\frac{d\sigma^{\pm}(\eta)}{||x-\eta||^{d-1}},\, x\in rS^d.
\]
For $x\in rS^d$, write $x=r\zeta,\, \zeta\in S^d$ and let $\eta\in S^d$. Then by (2.3) and (2.5)
\begin{eqnarray}\label{eq 3.2}
&& \frac{1}{||x-\eta||^{d-1}}=\frac{1}{(1+r^2-2r(\zeta.\eta))^{(d-1)/2}} \\
\nonumber && =\sum_{l}r^lP_l^{(d-1)/2)}((\zeta.\eta)) \\
\nonumber && =\sum_{l}r^l{l+d-2\choose l}P_l(d+1,(\zeta.\eta)) \\
\nonumber && =\sum_{l,k}r^{l}
\frac{\omega_{d}}{Z(d,l)}{l+d-2\choose l}Y_{l,k}(\zeta)
Y_{l,k}(\eta) \\
\nonumber && =\sum_{l,k}r^{l}\frac{\omega_{d}}{Z(d,l)}{l+d-2\choose l}
Y_{l,k}(\eta)Y_{l,k}(\zeta).
\end{eqnarray}
Now by (2.1)
\[
{l+d-2\choose l}\frac{l+d-1}{2l+d-1}
\frac{l!(d-1)!}{(l+d-1)!}=\frac{d-1}{2l+d-1}.
\]
Thus we learn that
\begin{equation}\label{eq 3.3}
\frac{1}{||x-\eta||^{d-1}}=\sum_{l,k}\frac{(d-1)\omega_{d}}{2l+d-1}
r^{l}Y_{l,k}(\eta)Y_{l,k}(\zeta).
\end{equation}
For a given $l$ and $k$, let us set for convenience
\[
\hat{\sigma}^{\pm}(l,k):=\int_{S^d}Y_{l,k}(\eta)d\sigma^{\pm}(\eta),
\]
\[
\sigma_B^{\pm}(r\eta):=\sigma_{B,r}^{\pm}(\eta),\, \eta\in S^d
\]
and
\[
\hat{\sigma}^{\pm}_{B,r}(l,k):=\int_{S^d}Y_{l,k}(\eta)d\sigma_B^{\pm}(r\eta).
\]
Then we have
\begin{equation}\label{eq 3.4}
\int_{S^d}\frac{1}{||x-\eta||^{d-1}}d\sigma^{\pm}(\eta)
=\sum_{l,k}\frac{(d-1)\omega_{d}}
{2l+d-1}r^{l}\hat{\sigma}^{\pm}(l,k)Y_{l,k}(\zeta).
\end{equation}
Also, we claim that
\begin{eqnarray}\label{eq 3.5}
&&\int_{rS^d}\frac{d\sigma_B^{\pm}(y)}{||x-y||^{d-1}}=\int_{S^d}
\frac{d\sigma_B^{\pm}(r\eta)}{||r\zeta-r\eta||^{d-1}} \\
\nonumber && =\sum_{l,k}\frac{(d-1)\omega_d}{(2l+d-1)r^{d-1}}
\hat{\sigma}_{B,r}^{\pm}(l,k)Y_{l,k}(\zeta).
\end{eqnarray}
To see (3.5), we proceed as follows. Consider the Newtonial
kernel on $S^d$ defined by way of the function
\[
g(t):=2^{(-d+1)/2}(1-t)^{(-d+1)/2},\, t\in [-1,1]
\]
and define for each $\delta>0$, a sequence of functions $g(t-\delta),\, t\in [-1,1]$.
This sequence coverges monotonically to $g(t)$ as $\delta\to 0^+$ for each
fixed $t$. It is known, that the
coefficients in the
expansion of $g(t-\delta)$ (with respect to $P_l(t,d+1)$) converge to
${l+d-2\choose l}$ as $\delta\to 0^+$.
Thus using the finiteness of the integral
\[
\int_{S^d}\frac{d\sigma^{\pm}(\eta)}{||z-\eta||^{d-1}},\, z\not\in S^d
\]
and the monontone convergence theorem together with the calculations of (3.3) yield
(3.5). Comparing coefficients in the two formulas (3.4) and (3.5) above, we learn that
\begin{equation}\label{eq 3.6}
\hat{\sigma}_{B,r}^{\pm}(l,k)=r^{l+d-1}\hat{\sigma}^{\pm}(l,k),\,
k=1,...,Z(d,l),\, l\geq 0.
\end{equation}
Armed with (3.6), we now recall (see (2.16)), that
\[
D(f_r)(y)=\sum_{l,k}\frac{(2l+d-1)}{(d-1)\omega_d}\hat{f_r}(l,k)Y_{l,k}(y),\, y\in S^d.
\]
Then we have on the one hand, by (3.6) and the property of balayage,
\begin{eqnarray}\label{eq 3.7}
&& \int_{S^d}f(y)d\sigma^{\pm}(y)=\int_{rS^d}f(y)d\sigma_B^{\pm}(y) \\
\nonumber && =\int_{S^d}f(r\eta)d\sigma_B^{\pm}(r\eta)=\int_{S^d}
f_r(\eta)d\sigma_{B,r}^{\pm}(\eta) \\
\nonumber && =\sum_{l,k}\hat{f_r}(l,k)\int_{S^d}Y_{l,k}(\eta)
d\sigma_{B,r}^{\pm}(\eta) \\
\nonumber && =\sum_{l,k}\hat{f_r}(l,k)\hat{\sigma}_{B,r}^{\pm}(l,k) \\
\nonumber && =\sum_{l,k}r^{l+d-1}\hat{f_r}(l,k)\hat{\sigma}^{\pm}(l,k).
\end{eqnarray}
On the other hand, using orthogonality and the definition of $D$, we
see that
\begin{eqnarray}\label{eq 3.8}
&& \frac{1}{r}\int_{rS^d}D(f)(y)U^{\sigma^{\pm}}(y)d\mu_d(y)=
\frac{1}{r}\int_{S^d}D(f)(r\eta)U^{\sigma^{\pm}}(r\eta)
r^dd\mu_d(\eta) \\
\nonumber && =\frac{1}{r}\int_{S^d}D(f_r)(\eta)U^{\sigma^{\pm}}(r\eta)r^dd\mu_d(\eta) \\
\nonumber && =\frac{1}{r}\int_{S^d}
\left(\sum_{l,k}\hat{f_r}(l,k)\frac{(2l+d-1)}{(d-1)\omega_d}Y_{l,k}(\eta)\right)
\left(\sum_{m,j}\frac{(d-1)\omega_d}{2m+d-1}r^{m}\hat{\sigma}^{\pm}(m,j)
Y_{m,j}(\eta)r^d\right)d\mu_d(\eta) \\
\nonumber && =\sum_{l,k}r^{l+d-1}\hat{f_r}(l,k)\hat{\sigma}^{\pm}(l,k).
\end{eqnarray}
Comparing (3.7) and (3.8) gives Lemma 1. $\Box$
\medskip

We remark that the proof of Lemma 1 shows that the expansion of $f_r$ needed
in the proof is only needed $\sigma^{\pm}_{B,r}\, a.e$.
\medskip

{\bf The Proof of Theorem 2}\, Theorem 2 follows from Lemma 1 and Minkowski's
inequality. $\Box$
\medskip

We now proceed with the
\medskip

{\bf Proof of Proposition 1}\, We first establish (2.17). We may assume with loss of generality that $p=\infty$. Using the addition theorem,
the Cauchy-Swartz inequality, and the estimate,
\[
Z(d,l)\leq e^{d}l^{d-1},\, l\geq 0
\]
we see that we have for $x\in S^d$,
\begin{eqnarray*}
&& |f_r(x)|^2=\left|\sum_{l,k}\hat{f_r}(l,k)Y_{l,k}(x)\right|^2\\
&& \leq \left(\sum_{l,k}(\hat{f_r}(l,k))^2m_l^2\right)
\left(\sum_{l,k}Y_{l,k}^2m_l^{-2}\right)\\
&& \leq ||f_r||_{H_s(S^d)}^{2}\left[\sum_{l=1}^{\infty}\frac{Z(l,d)m_l^{-2}}{\omega_d}+
\frac{1}{\omega_d}\right] \\
&& \leq ||f_r||_{H_s(S^d)}^{2}
\left[\sum_{l=1}^{\infty}\frac{e^d\omega_{d}(d-1)^2l^{d-1-2s}}{(2l+d-1)^2}+
\frac{1}{\omega_d}\right].
\end{eqnarray*}
Note that the last sum is finite since $s>\frac{d}{2}$. Thus (2.17) holds. To see (2.18), we proceed
much as in (2.17). As $s>\frac{d}{2}$, we readily have for $x\in S^d$
\begin{eqnarray*}
&& \left|\sum_{l,k}\frac{(2l+d-1)}{(d-1)\omega_d}\hat{f_r}(l,k)Y_{l,k}(x)\right|^2 \\
&& \leq \left(\sum_{l,k}(\hat{f_r}(l,k))^2m_l^2\right)
\left(\sum_{l,k}\frac{(2l+d-1)^2}{(d-1)^2(\omega_d)^2}Y_{l,k}^2m_l^{-2}\right)\\
&& \leq ||f_r||_{H_s(S^d)}^{2}
\left[\sum_{l=1}^{\infty}\frac{e^dl^{d-1-2s}}{\omega_d}+\frac{1}{\omega_d^3}\right].
\end{eqnarray*}
It remains to show (2.19).
In light of (2.17) and (2.18), it is easy to see that we have for any
$\zeta$ and $\eta$
on $S^d$ the estimate:
\begin{eqnarray}\label{eq 3.9}
&& |f_r(\eta)-f_r(\zeta)|\leq C||f_r||_{H_s(S^d)}\times \\
\nonumber && \times \left[\sum_{l,k}\frac{(d-1)^2\omega_d^2}{(2l+d-1)^2}
(Y_{l,k}(\eta)-Y_{l,k}(\zeta))^2l^{-2s}+
\sum_{k}(Y_{0,k}(\eta)-Y_{0,k}(\zeta))^2\right]^{1/2}.
\end{eqnarray}
Now it is a straightforward consequence of the addition formula, that for
every fixed $l\geq 0$, we have
\begin{eqnarray*}
&& \sum_{k=1}^{Z(d,l)}(Y_{l,k}(\zeta)-Y_{l,k}(\eta))^2 \\
&& =\frac{2Z(d,l)}{\omega_d}\left(1-P_l(\eta.\zeta,d+1)\right).
\end{eqnarray*}
Thus inserting the above estimate into (3.9) yields:
\begin{eqnarray}\label{eq 3.10}
&& |f_r(\eta)-f_r(\zeta)|\leq C_1||f_r||_{H_s(S^d)} \times \\
\nonumber && \times
\left[\sum_{l=1}^{\infty}l^{-2s-3+d}\left|1-P_l(\eta.\zeta,d+1)\right|
+ C_2\left|1-P_0(\eta.\zeta,d+1)\right|\right]^{1/2}.
\end{eqnarray}
We need to estimate the righthand side of (3.10). First,
we recall, (see (2.3), \cite[pp 63]{Sz} and \cite[pp 170]{Sz}),
that if $P_l^{(d/2-1,d/2-1)}$ is the Jacobi polynomial of degree $l\geq 0$, then
we have for every $l\geq 0$ and $x\in [-1,1]$,
\[
P_l(x,d+1)={l+d/2-1\choose l}^{-1}P_l^{(d/2-1,d/2-1)}(x)
\]
and
\[
(P_l^{(d/2-1,d/2-1)})'(1) \leq
\left \{ {
\begin{array}{ll}
Cl^{1+d/2}, & l>0 \\
C_1, & l=0
\end{array} }
\right.
\]
Thus we deduce that for every $l\geq 1$,
\begin{equation}\label{eq 3.11}
{\rm max}_{|t|\leq 1}P_l'(t,d+1)=P_l'(1,d+1)\leq
\left \{ {
\begin{array}{ll}
Cl^2,& l>0 \\
C_1, & l=0
\end{array} }
\right.
\end{equation}
We also have by the first mean value theorem and our normalization, that
for every $l\geq 0$,
\begin{equation}\label{eq 3.12}
1-P_l(\eta.\zeta,d+1)=\frac{1}{2}P'_l(t,d+1)2(1-\eta.\zeta)=\frac{1}{2}P'_l(t,d+1)||\eta-\zeta||^2
\end{equation}
for some $t\in [-1,1]$. Thus, we have
\[
\sum_{l=1}^{\infty}l^{-2s-3+d}|1-P_l(\eta.\zeta,d+1)|\leq C\left\{\sum_{l=1}^{\infty}l^{-2s-3+d}\right\}||\eta-\zeta||^2.
\]
(2.19) now follows. $\Box$
\bigskip

{\bf The Proof of Theorem 3}\, This follows in a straightforward way from Theorem 2.
$\Box$

\setcounter{equation}{0}
\section{Proof of Theorem 4}

In this last section, we establish Theorem 4.
\medskip

{\bf Proof of Theorem 4}\, Let us write $\sigma_n=\sigma:=\mu-\nu_n$ and take
$x\in S^d$. First we observe that we have
\begin{eqnarray}\label{eq 4.1}
&& -U^{-\nu_n}(rx)=\int_{S^d}\frac{1}{||rx-y||^{d-1}}d\nu_n(y) \\
\nonumber && =\int_{S^d}\frac{1}{||rx-y||^{d-1}}\sum_{k=1}^na_{k,n}
d\delta_{t_{k,n}}(y)=
\sum_{k=1}^n\frac{a_{k,n}}{||rx-t_k||^{d-1}} \\
\nonumber && =\sum_{k=1}^n\frac{\mu(R_{k,n})}
{||rx-t_{k,n}||^{d-1}}=\sum_{k=1}^n\int_{R_{k,n}}
\frac{1}{||rx-t_{k,n}||^{d-1}}d\mu(y).
\end{eqnarray}
On the other hand,
\begin{equation}\label{eq 4.2}
U^{\mu}(rx)=\int_{S^d}\frac{d\mu(y)}{||rx-y||^{d-1}}=\sum_{k=1}^n\int_{R_{k,n}}
\frac{d\mu(y)}{||rx-y||^{d-1}}.
\end{equation}
Thus, using (4.1) and (4.2), we have
\begin{equation}\label{eq 4.3}
\left|U^{\sigma}(rx)\right| \leq \sum_{k=1}^n\int_{R_{k,n}}\left|
\frac{1}{||rx-y||^{d-1}}-\frac{1}{||rx-t_k||^{d-1}}\right|d\mu(y).
\end{equation}
Now let us estimate the integrand in (4.3). We first write for $x,y\in S^d$:
\[
||rx-y||^{d-1}=(r^2-2r(x.y)+1)^{(d-1)/2}.
\]
A simple calculation also shows that
\begin{eqnarray*}
&& \frac{d}{d\theta}(1+r^2-2r\theta)^{-(d-1)/2} \\
&& =(d-1)r((1+r^2-2r\theta)^{-(d+1)/2} \\
&& \leq (d-1)(1-r)^{-(d+1)}
\end{eqnarray*}
for any $-1\leq \theta\leq 1$. Thus an application of the mean value theorem to
(4.3) easily yields
\begin{eqnarray*}
&& \left|U^{\sigma}(rx)\right|\leq (d-1)||\mu||(1-r)^{-(d+1)}||y-t_k|| \\
&& \leq (d-1)||\mu||(1-r)^{-(d+1)}||{\cal R}||.
\end{eqnarray*}
Theorem 4(a) is then proved. Next, it is known that, we can find a reduction
$E$ of $E_o$ and
an explicit finite disjoint partitioning
${\cal R}$ of $S^d$ with each $R\in {\cal R}$ containing one point of $E$ and
at least one
point of $E_o$. Moreover,
\begin{equation}\label{eq 4.4}
\delta_{E_o}\leq \delta_{E}<||{\cal R}||\leq 8d\sqrt{2d(d+1)}\delta_{E_o}.
\end{equation}
Then (2.27) follows from (2.26). Moreover, (2.28) then follows from
(2.25) and (2.27). This completes the proof of Theorem 4. $\Box$
\bigskip

{\bf Acknowlegement}\, Support from the American Mathematical Society is acknowledged.

\bigskip

\end{document}